\documentclass[10pt,twoside]{article}

\usepackage{amssymb,color,ulem}

\usepackage[usenames,dvipsnames]{xcolor}

\usepackage{graphicx}

\def\qed{\hfill $\square$}

\usepackage{amsthm,verbatim}
\theoremstyle{definition}

\newtheorem{theorem}{Theorem}[section]
\newtheorem{cor}{Corollary}[section]

\newtheorem{ex}{Example}[section]
\newcommand{\eps}{\varepsilon}

\usepackage{anysize}  %----------------margins
\usepackage{bigstrut} % controling row heights of tables
\usepackage{mathrsfs}
 \usepackage{amsmath,amssymb}

\makeatletter
\evensidemargin \oddsidemargin
\makeatother

\begin{document}
\thispagestyle{empty}
\setcounter{page}{1}

\noindent
{\footnotesize {\rm To appear in\\[-1.00mm]
{\em Dynamics of Continuous, Discrete and Impulsive Systems}}\\[-1.00mm]
http:monotone.uwaterloo.ca/$\sim$journal} $~$ \\ [.3in]

%%USE THE STUFF BELOW AS A GUIDE TO SET UP THE START OF PAPER%%

\begin{center}
{\large\bf Dwell time for switched systems with multiple equilibria on a finite
time-interval}

\vskip.20in

Oleg Makarenkov \ and \ Anthony  Phung \\[2mm]
{\footnotesize
Department of Mathematical Sciences\\
University of Texas at Dallas, Richardson, TX 75080, USA \\
}
\end{center}

{\footnotesize
\noindent
{\bf Abstract.} We describe the behavior of solutions of switched systems with multiple globally exponentially stable equilibria. We introduce an ideal attractor and show that the solutions of the switched system stay in any given $\eps$-inflation of the ideal attractor if the frequency of switchings is slower than a suitable dwell time $T$. In addition, we give conditions to ensure that the $\eps$-inflation is a global attractor. Finally, we investigate the effect of the increase of the number of switchings on the total time that the solutions need to go from one region to another.
 \\[3pt]
}

\vskip.2in

\section{Introduction}

Dwell time is the lower bound on the time between successive switchings of the switched system
\begin{equation}\label{ss}
   \dot x=f_{u(t)}(x),\qquad u(t)\mbox{ is a piecewise constant function, } x\in\mathbb{R}^n,
\end{equation}
which ensures a required dynamic behavior under the assumption that each of the subsystems 
\begin{equation}\label{ssi}
   \dot x=f_{u}(x),\qquad u\in\mathbb{R},\ 
x\in\mathbb{R}^n,
\end{equation}
 possess a globally stable equilibrium $x_u.$ When all the equilibria $\{x_{u(t)}\}_{t\ge t_0}$ coincide, the dwell time $T>0$ which gives global exponential stability of the common equilibrium $x_0$ is computed e.g. in Liberzon \cite[\S3.2.1]{liberzon}. Specifically, the result of \cite[\S3.2.1]{liberzon} gives a formula for $T$ which makes $x_0$ globally exponentially stable for any piecewise constant function $u(t)$ whose discontinuities $t_1,t_2,\ldots$ verify 
\begin{equation}\label{p}
  |t_i-t_{i-1}|\ge T.
\end{equation}
The case where the equilibria are distinct is covered in Alpcan-Basar \cite{alpcan}, who offered a dwell time $T$ that ensures global exponential stability of a suitable set $A\supset \{x_{u(t)}\}_{t\ge t_0}$ for any $u(t)$ whose discontinuities verify (\ref{p}). The problem of stability of switched systems with multiple equilibria appears e.g. in differential games, load balancing, agreement and robotic navigation (see \cite{alpcan,spong} and references therein).

\vskip0.2cm

\noindent A deeper analysis of the dynamics of switched systems with multiple equilibria was recently carried out in Xu et al \cite{xu}, who gave a sharp formula for the attractor $A$ in the case of quasi-linear switched systems (\ref{ss}). Assuming that $u(t)$ is periodic and denoting by $t\mapsto X_u(t,x)$ the solution of (\ref{ssi}) with the initial condition $X_u(0,x)=x$, the paper \cite{xu} investigated the asymptotic attractivity of 
$$
  A=\bigcup_{t\ge t_0,\ \tau\ge t_0}\{X_{u(\tau)}(t,x_{u(\tau)})\}.
$$
\noindent The motivation for our paper comes from the problem of planning the motion of a 3-D walking robot, where "turn left", "walk straight" and "turn right" correspond to $u(t)=-1$, $u(t)=0$ and $u(t)=1$ respectively, see Gregg et al \cite{gregg}. It is not the asymptotic attractivity of $A$ which is of importance for the robot turning maneuver but rather an appropriate attractivity of $A$ during the time of the maneuver.  The goal of this paper is to provide a dwell time which can ensure the required attractivity.

\vskip0.2cm

\noindent The paper is organized as follows. In the next section of the paper, we prove our main result (Theorem~\ref{thm1}). Given $\eps>0$, Theorem~\ref{thm1} provides a dwell time $T>0$ such that the solutions of (\ref{ss}) with the initial conditions in the $\eps$-neighborhood $B_\eps(A)$ of $A$ never leave $B_\eps(A)$ in the future. Theorem~\ref{thm1} can be viewed as a version of \cite[Theorem~1]{xu} for fully nonlinear systems. In section 3, we compute (Theorem~\ref{thm2}) a dwell time to ensure that the attractor $B_\eps(A)$ is reached asymptotically from any initial condition. 
The proof of Theorem~\ref{thm2} follows the ideas of Alpcan-Basar \cite{alpcan}. However, we offer weaker conditions where the Lyapunov functions of subsystems (\ref{ssi}) are not supposed to respect any uniform estimates. A particular case study where the Lyapunov functions of subsystems (\ref{ssi})  are shifts of one another is addressed in section~4. In this section, we consider a switched system which switches between two subsystems $u=u_1$ and $u=u_2$ and analyze the solutions of the switched system with the initial conditions in $B_\eps(A)$. Let $x_1$ and $x_2$ be the equilibria of subsystems $u=u_1$ and $u=u_2$ respectively.
  The result of section~4 (Theorem~\ref{thm41}) clarifies whether or not the solutions from the neighborhood of $x_1$ reach the neighborhood of $x_2$ faster if the switching signal is amended in such a way that an additional switching occurs between $u=u_1$ and $u=u_2.$ In  other words, section~4 investigates whether or not adding more discrete events is alone capable of making the dynamics inside $B_\eps(A)$ faster. Examples~\ref{ex1} and \ref{ex2} illustrate the conclusions of Theorems~\ref{thm1} and \ref{thm41}.

\section{The local trapping region}

Let $x_u$ be the unique equilibrium of (\ref{ssi}). We assume that for any $u$, system (\ref{ssi}) admits a global Lyapunov function $V_u$ such that
\begin{eqnarray}
&&   \alpha_u(\|x-x_u\|)\le V_u(x)\le \beta_u(\|x-x_u\|),\qquad x\in\mathbb{R}^n,\label{ab}\\
&& (V_u)'(x)f_u(x)\le -k_u V_u(x),\qquad x\in\mathbb{R}^n,\label{k}
\end{eqnarray}
where $\alpha$, $\beta$ are strictly monotonically increasing functions with $\alpha_u(0)=\beta_u(0)$, and $k_u>0.$
Introduce the following trapping regions
\begin{equation}\label{Ni}
\renewcommand*{\arraystretch}{1.5}
\begin{array}{rcl}
  N_u^\eps&=&   \left\{x:V_u(x)\le \eps \right\},\\
  L_{u_1,u_2}^\eps (t)& = & \bigcup\limits_{x\in N_{u_1}^\eps}\left\{X_{u_2}(t,x)\right\}
\end{array}
\end{equation}
and define the dwell time that the solutions need to go from $N_{u_1}^\eps$ to $N_{u_2}^\eps$ as
%Given a function $u:[t_0,\infty)\to \overline{1,m}$ and the respective increasing sequence $(t_1,t_2,\ldots)=\{t_i\}_{i\in I}$, let $u_i=u(t_i)$ and 
\begin{equation}\label{Ti}
T_{u_1,u_2}^\eps=-\dfrac{1}{k_{u_2}}\ln\dfrac{\eps}{\beta_{u_2}\left(\|x_{u_2}-x_{{u_1}}\|+\alpha_{{u_1}}^{-1}(\eps)\right)}.
\end{equation}

\begin{theorem}\label{thm1} Assume that 
\begin{itemize}
\item[(A1)] $f_u\in C^1(\mathbb{R}^n,\mathbb{R}^n)$ for any $u\in\mathbb{R}$,
\item[(A2)] for any $u\in\mathbb{R}$,  system (\ref{ssi}) admits an equilibrium $x_u$ whose Lyapunov function $V_u$ satisfies (\ref{ab})-(\ref{k}), where $\alpha,\beta\in C^0(\mathbb{R},\mathbb{R})$ are strictly increasing functions, $\alpha(0)=\beta(0)=0,$ $k_u>0,$
\item[(A3)] $u:[t_0,\infty)\to \mathbb{R}$ is a piecewise constant function.
\end{itemize}
Let $\{t_1,t_2,\ldots\}=\{t_i\}_{i\in I}$ be a finite or infinite increasing sequence of points of discontinuity of $u$ and 
$$
    u_i=u(t_i+0).
$$
If
$$
   t_i-t_{i-1}\ge T_{u_{i-1},u_i}^\eps,
$$
then for any solution $x$ of (\ref{ss}) with
$$
   x(t_{i-1})\in N_{u_{i-1}}^\eps,\quad i\in I,
$$
one has
\begin{eqnarray}
& x(t)\in L_{u_{i-1},u_i}^\eps(t),  &\ \ t_{i-1}\le t\le t_i,\ \ \ i\in I,\label{prove1} \\
&x(t_i)\in N_{u_i}^\eps, & \ \ i\in I. \label{prove2}
\end{eqnarray}
\end{theorem}

\noindent {\bf Proof.} We only have to prove (\ref{prove2}) because the validity of (\ref{prove1}) follows directly from the definition of $L_{u_{i-1},u_i}^\eps$. Let $x(t_{i-1})\in N_{u_{i-1}}^\eps.$ Our goal is to show that $x(t_i)\in N_{u_i}^\eps.$
Given $\eps>0$, define $\delta>0$ as
$$
   \delta=\beta_{u_i}\left(\|x_{u_i}-x_{u_{i-1}}\|+\alpha_{u_{i-1}}^{-1}(\eps)\right).
$$ 
\begin{figure}[h]\center
\includegraphics[scale=0.7]{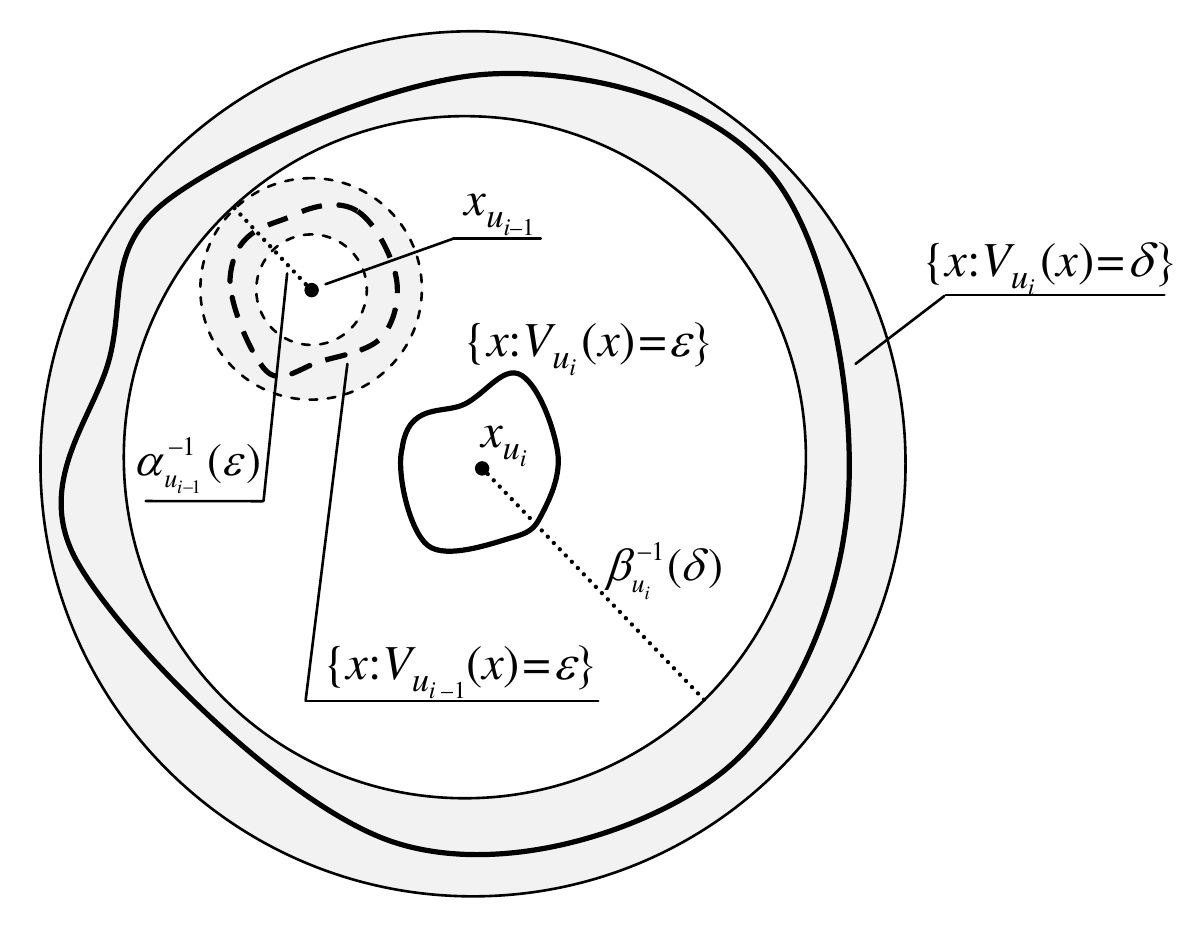}\ \ \ \
\vskip-0.2cm
\caption{\footnotesize Illustration of the proof of Theorem~\ref{thm1}. The two gray rings are the estimates for the sets $\{x:V_{u_{i-1}}(x)=\eps\}$ and $\{x:V_{u_{i}}(x)=\delta\}$ given by condition (\ref{ab}).} \label{fig11}
\end{figure}
\noindent By construction (see Fig.~\ref{fig11}), 
$N_{u_i}^\delta\supset N_{u_{i-1}}^\eps$, and so 
$x(t_{i-1})\in N_{u_i}^\delta$. 
%it is sufficient to show that for any solution $x$ of (\ref{ssi}) with  $u=u_i$ and $x(0)\in N_{u_i}^\delta$ it takes at most $T_i^\eps$ to get to $N_{u_i}^\eps.$ 
 Introduce
$$
   v(t)=V_{u_i}(x(t)).
$$
By (\ref{k}) we have
\begin{eqnarray*}
  \dot v(t)&\le &-k_i v(t),\quad t_{i-1}\le t\le t_i,\\
   v(t_{i-1})&\le &\delta.
\end{eqnarray*}
By the comparison lemma (see e.g. \cite[Lemma~16.4]{amann}), it holds that
$$
   v(t)\le p(t),
$$
where $p(t)$ is the solution of
\begin{eqnarray*}
  \dot p(t)&= &-k_i p(t),\quad t_{i-1}\le t\le t_i,\\
   p(t_{i-1})&= &\delta.
\end{eqnarray*}
At the same time,
$$
  p(t_i)=e^{-{k_{u_i}}(t_i-t_{i-1})}\delta \le e^{-k_{u_i} T_{u_{i-1},u_i}^\eps}\delta=\dfrac{\eps}{\beta_{u_i}\left(\|x_{u_i}-x_{u_{i-1}}\|+\alpha_{u_{i-1}}^{-1}(\eps)\right)}\delta=\eps.
$$
Therefore, $V_{u_i}(x(t_i))\le\eps$, which completes the proof.\qed

\vskip0.2cm

\noindent Theorem~\ref{thm1} suggests the following definition of the $\eps$-inflation $A_\eps$ of the ideal attractor $A$ of (\ref{ss}). Given a function $u:[t_0,\infty)\to \mathbb{R}$ and the respective increasing sequence $(t_1,t_2,\ldots)=\{t_i\}_{i\in I}$, let $u_i=u(t_i)$ and 
$$
\renewcommand*{\arraystretch}{1.5}
   A_\eps(t)=\left\{\begin{array}{lll}
L_{u_{i-1},u_i}^\eps(t), & t_{i-1}\le t<t_i, & i\in I,\\
L_{u_{\max (I)-1},u_{\max (I)}}^\eps(t),
& t\ge t_{\max (I)}, & {\rm if}\ I\ {\rm is\ finite.}\end{array}\right.
$$
\begin{cor}\label{Tloc} Let the assumptions (A1)-(A3) of Theorem~\ref{thm1} hold.  If 
$$
   t_i-t_{i-1}\ge  
\sup\limits_{i\in I} T_{u_{i-1},u_i}^\eps=:T_{loc}^\eps, \quad i\in I,
$$
then, for any solution $x$ of (\ref{ss}) with the initial condition
$$
   x(t_{0})\in A_\eps (t_{0}),
$$ 
 one has
$$
   x(t)\in A_\eps(t),\quad t\ge t_{0}.
$$
\end{cor}
\noindent Note, $\sup\limits_{i\in I} T_{u_{i-1},u_i}^\eps$ is finite when $t\mapsto u(t)$ takes a finite number of values on $[t_0,\infty).$ 

\begin{ex}\label{ex1} To illustrate Theorem~\ref{thm1}, we consider the following switched system (slightly modified from Example~2 in \cite{alpcan})
\begin{equation}\label{ss_ex1}
\dot x=\left(\begin{array}{cc}
  -1 & -1 \\
 1 & -1 \end{array}\right)x+\left(\begin{array}{c} u \\ 1\end{array}\right),
\end{equation}
whose unique equilibrium is given by
$$
   x_u=\dfrac{1}{2}\left(\begin{array}{c}
u-1\\ u+1
\end{array}
\right).
$$
Introduce the three discrete states $u_1,$ $u_2,$ and $u_3$ as
$$
   u_1=1,\qquad u_2=0,\qquad u_3=-1,
$$
and consider 
$$
  \eps=0.05.
$$
If the Lyapunov function $V_u(x)$ is selected as
$$
   V_u(x)=\|x-x_u\|^2,
$$
then formulas (\ref{Ni}) and (\ref{Ti}) yield 
$$
   N_u^\eps=\left\{x:\|x-x_u\|\le \sqrt{\eps}\right\},\quad T_{loc}^\eps\approx 1.426.
$$
Therefore, for the control input 
\begin{equation}\label{ci}
   u(t)=\left\{\begin{array}{l}u_1,\quad t\in[0,T),\\
u_2,\quad t\in[T,2T),\\
u_3,\quad t\ge 2T,\end{array}\right.\qquad T=1.43,
\end{equation}
and for any solution $x$ of (\ref{ss_ex1}), 
Theorem~\ref{thm1} ensures the following: 
$$
   {\rm if} \ x(0)\in N^\eps_{u_1}, \ {\rm then} \ x(T)\in N^\eps_{u_2} \ {\rm and} \ x(2T)\in N^\eps_{u_3}.
$$
\end{ex}

\begin{figure}[h]\center
\includegraphics[scale=0.44]{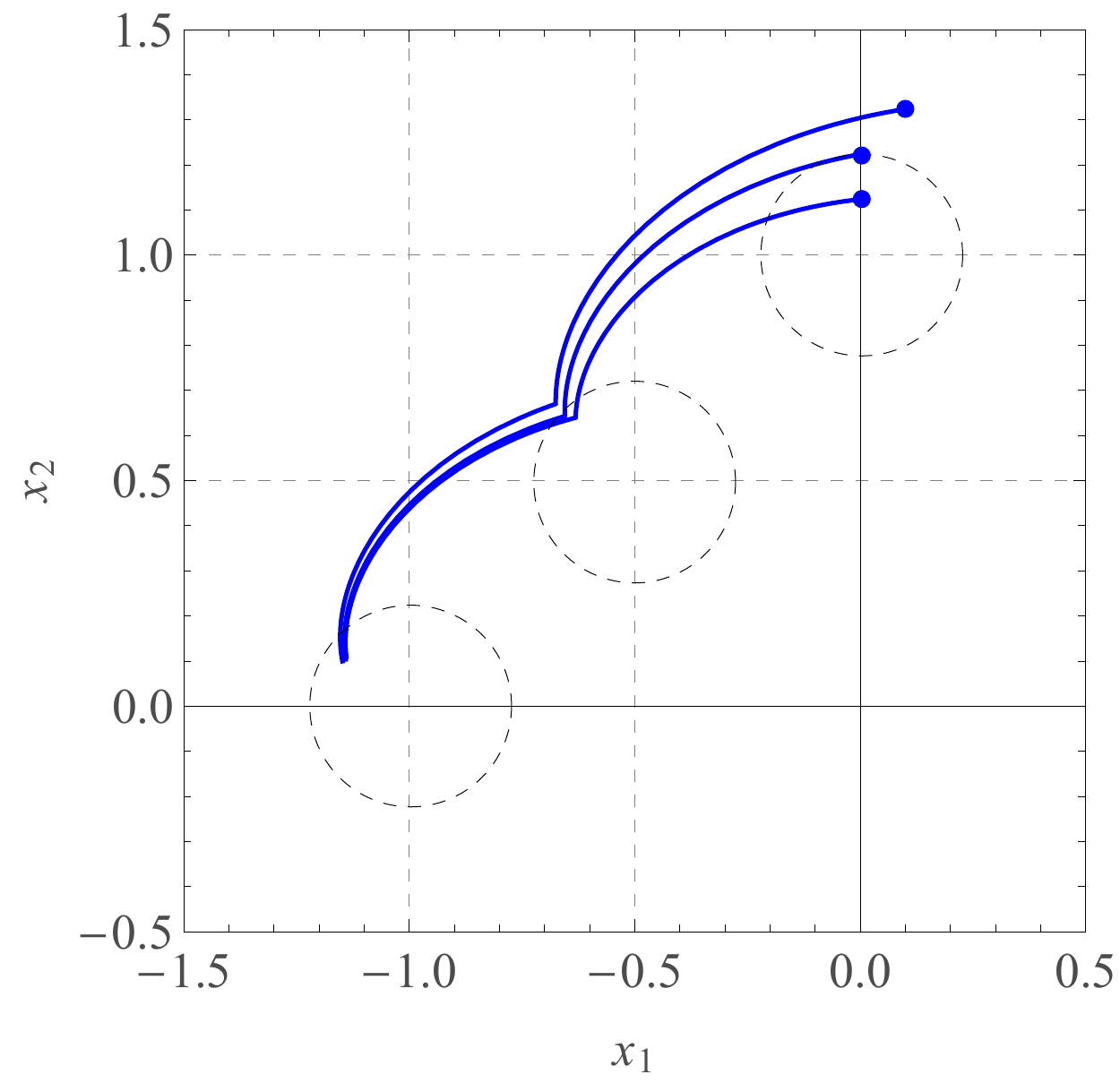}\ \ \ \
\includegraphics[scale=0.44]{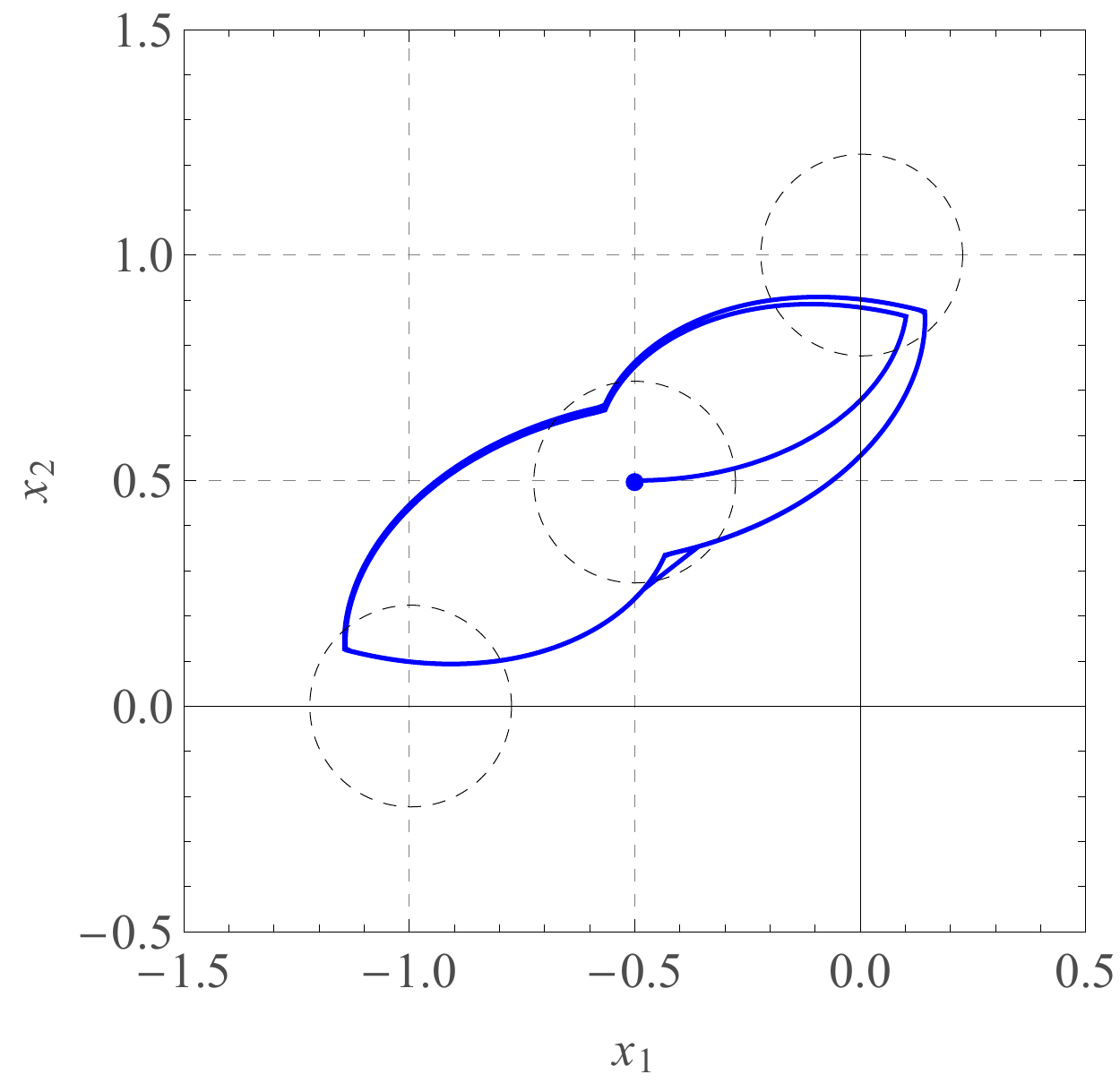}
\vskip-0.4cm
\caption{\footnotesize Left: Solutions of switched system (\ref{ss}) with initial conditions (blue dots) inside $N^{0.05}_{u_1}$, on the boundary of $N^{0.05}_{u_1}$, and outside $N^{0.05}_{u_1}$, and for the control input $u(t)$ given by (\ref{ci}). Right: The solution of switched system (\ref{ss}) with the initial condition $(-0.5,0.5)^T$ for the $4T$-periodic control input (\ref{4Tci}).} \label{fig1}
\end{figure}

\noindent Figure~\ref{fig1}(left) documents the sharpness of  the dwell time $T$. Indeed, the figure shows that if the initial condition $x(0)$ deviates to the outside of  $N^\eps_{u_1}$ just a little bit, then the dwell time $T$ is no longer sufficient to get $x(T)\in N^\eps_{u_2}$ (though we still have $x(2T)\in N^\eps_{u_3}$ for this solution).

\vskip0.2cm

\noindent To demonstrate that trapping regions $N^\eps_{u_1}$, $N^\eps_{u_2}$, $N^\eps_{u_3}$ (and thus the $\eps$-inflated attractor $A_\eps$, see Corollary~\ref{Tloc}) provide a rather sharp estimate for the location of the attractor of (\ref{ss_ex1}), we extend the input $u(t)$ to $[0,4T]$ as 
\begin{equation}\label{4Tci}
   u(t)=\left\{\begin{array}{l}u_1,\quad t\in[0,T),\\
u_2,\quad t\in[T,2T),\\
u_3,\quad t\in[2T,3T),\\
u_2,\quad t\in[3T,4T),
\end{array}\right.
\end{equation}
and then continue it to the entire $[0,\infty)$
by $4T$-periodicity. The respective solution $x$ of (\ref{ss_ex1}) with the initial condition $x(0)=(-0.5,0.5)^T$ is plotted in Fig.~\ref{fig1} (right). The drawing shows that the switching points of the solution $x$ are very close to the boundaries of the trapping regions $N^\eps_{u_1}$, $N^\eps_{u_2}$, $N^\eps_{u_3}$, i.e. there is only a  little window to reduce the size of those regions.

\section{Global attractivity of the local trapping region} 

\begin{theorem} \label{thm2} 
  Let the assumptions (A1)-(A3) of Theorem 2.1 hold and $I$ be infinite.
  Fix $\varepsilon>0$ and suppose that there exists constants $\mu_i(\varepsilon)$ such that
  \begin{equation*}
  \frac{V_{u_{i+1}}(x)}{V_{u_i}(x)}\leq \mu_i(\varepsilon), \quad x\in\mathbb{R}^n\setminus N_{u_i}^{\varepsilon}, \quad i\in\mathbb{N}\cup\{0\}.
\end{equation*}
Finally, assume that
\begin{equation*}
  \mu_0(\varepsilon)\cdot\dotsc\cdot\mu_i(\varepsilon)e^{\displaystyle -\int_{t_0}^{t_{i+1}}k_{u(s)}ds}\to 0, \quad \text{as } i\to\infty.
\end{equation*}
Then, $x(\hat{T})\in N_{u_i}^{\varepsilon}$ for some $\hat{T}>0$ and some $i\in\mathbb{N}$.
\end{theorem}

\noindent{\bf Proof.}
  Let $W(t)=e^{k_{u(t)}t}V_{u(t)}(x(t))$, where $t\in[t_0,\infty)$.
  Then, for $t\in[t_i,t_{i+1})$,
  \begin{equation*}
    W'(t)=k_{u_i}W(t)+e^{k_{u_i}t}\frac{d}{dt}V_{u_i}(x(t))\leq k_{u_i}W(t)-k_{u_i}e^{k_{u_i}t}V_{u_i}(x(t))=0,
  \end{equation*}
  which means that $W$ is decreasing on $[t_i,t_{i+1})$.
  In particular,
  \begin{equation*}
    W(t_i^+)\geq W(t_{i+1}^-).
  \end{equation*}
  On the other hand,
  \begin{equation*}
    \frac{W(t_{i+1}^+)}
    {W(t_{i+1}^-)}
    =
    \frac{e^{k_{u_{i+1}}t_{i+1}}}
    {e^{k_{u_i}t_{i+1}}}
    \cdot
    \frac{V_{u_{i+1}}(x(t_{i+1}))}
    {V_{u_i}(x(t_{i+1}))}
    \leq
    \frac{e^{k_{u_{i+1}}t_{i+1}}}
    {e^{k_{u_i}t_{i+1}}}
    \cdot
    \mu_i(\varepsilon)=:\tilde{\mu}_i.
  \end{equation*}
  Therefore,
  \begin{equation}\label{eq:tildeInequality1}
    \tilde{\mu}_iW(t_i^+)
    \geq
    \tilde{\mu}_iW(t_{i+1}^-)
    \geq
    W(t_{i+1}^+)
  \end{equation}
  Replacing $i$ by $i-1$ and combining with (\ref{eq:tildeInequality1}), one gets $\tilde{\mu}_{i-1}\tilde{\mu}_i W(t_{i-1}^+)\geq W(t_{i+1}^+)$.
  Continuing this process for $i-2,i-3,$ etc. we obtain
  \begin{equation*}
    \tilde{\mu}_0\cdot\dotsc\cdot\tilde{\mu}_iW(t_0^+)
    \geq
    W(t_{i+1}^+).
  \end{equation*}
  Applying $e^{-k_{u_{i+1}}t_{i+1}}$ yields
  \begin{equation*}
    \tilde{\mu}_0\cdot\dotsc\cdot \tilde{\mu}_i e^{-k_{u_{i+1}}t_{i+1}}W(t_0^+)
    \geq
    V_{u_{i+1}}(x(t_{{i+1}})),
  \end{equation*}
  or, equivalently,
  \begin{equation*}
    V_{u_0}(x(t_0))e^{-\sum_{j=0}^{i}k_{u_j}(t_{j+1}-t_j)}
    \mu_0(\varepsilon)\cdot\dotsc\cdot \mu_i(\varepsilon)
    \geq
    V_{u_{i+1}}(x(t_{i+1})).
  \end{equation*}
  The left-hand-side approaches 0 as $i\to\infty$ by the assumption of the theorem.
  Therefore, $V_{u_{i+1}}(x(t_{i+1}))\to0$ as $i\to\infty$.
  The proof is complete. \qed

\begin{cor} {\bf (Alpcan-Basar \cite{alpcan})} \label{alpcan} Let the assumptions (A1)-(A3) of Theorem~\ref{thm1} hold and $I$ be infinite.   Fix $\eps>0$ and suppose that there exists a constant $\mu(\eps)>1$ such that 
$$
   \dfrac{V_{u(t)}(x)}{V_{u(\tau)}(x)}\le\mu(\eps),\quad x\in\mathbb{R}^n\backslash N^\eps_{u(\tau)},\quad \tau,t\ge t_0.
$$
Finally, assume that  $k=\inf\limits_{t\ge t_0}k_{u(t)}>0$ and 
consider $T^\eps_{glob}$ satisfying
$$
   T^\eps_{glob}>\dfrac{\ln(\mu(\eps))}{k}.
$$
 If 
$$
   t_i-t_{i-1}\ge T^\eps_{glob},\quad i\in\mathbb{N},
$$
then, $x(\hat T)\in N_{u_i}^\eps$ for some $\hat T>0$ and some $i\in \mathbb{N}.$
\end{cor}

\noindent {\bf Proof.} Let $\gamma>0$ be such that $T^\eps_{glob}=\dfrac{\ln(\mu(\eps)+\gamma)}{k}$. Let $\mu_0(\eps),\ldots\mu_i(\eps)$ be as given by Theorem~\ref{thm2}. Then 
\begin{eqnarray*}   
&& e^{-\sum_{j=0}^i k_{u_{j}}(t_{j+1}-t_j)}\mu_0(\eps)\cdot\ldots\cdot\mu_i(\eps)\le e^{-kiT^\eps_{glob}}\mu(\eps)^i=\\
&& =\left(\mu(\eps)+\gamma\right)^{-i}\mu(\eps)^i=\left(\dfrac{\mu(\eps)}{\mu(\eps)+\gamma}\right)^i\to 0\ \ {\rm as}\ \ i\to\infty.\qed
\end{eqnarray*}

\begin{cor} Let the conditions of Corollary~\ref{alpcan} hold. Let $T^\eps_{loc}$ and $T^\eps_{glob}$ be those given by Corollaries 
\ref{Tloc} and \ref{alpcan}.
If
$$
   t_i-t_{i-1}\ge \max\left\{T^\eps_{loc},T^\eps_{glob}\right\},\quad i\in \mathbb{N},
$$
then, for any solution $x$ of (\ref{ss}), there exists $\hat T>t_0$ such that 
$$
   x(t)\in A_\eps(t),\quad t\ge \hat T.
$$
\end{cor}

\section{Dependence of the dwell time on the number of discrete states} 

\noindent Suppose that $u(t)$ switches from $u_0$ to $u_1$ at $t=t_0$. According to Theorem~\ref{thm1}, it takes at most time $T_{u_0,u_1}^\eps$ (see formula (\ref{Ti}))
for a trajectory $x$ of (\ref{ss}) to go from $N_{u_0}^\eps$ to $N_{u_1}^\eps.$ The next theorem shows that adding more discrete states between $u_0$ and $u_1$ makes the travel time from $N_{u_0}^\eps$ to $N_{u_1}^\eps$ longer.

\begin{theorem} \label{thm41}
Let the assumptions (A1)-(A2) of Theorem~\ref{thm1} hold and suppose $\alpha_u=:\alpha$, $\beta_u=:\beta$, $k_u=:k$ don't depend on $u$.  Fix $d>0$ and $r>0$. Then there exists $\eps_0>0$ such that
$$
   T_{u_0,u_1}^\eps<{T}_{u_0,v}^\eps+{T}_{v,u_1}^\eps, 
$$ 
for any 
%\begin{equation}\label{fu}
 %  u(t)=\left\{\begin{array}{ll}
 %    u_0, & t<t_0,\\
  %   u_1, & t\ge t_0,
  %   \end{array}\right.\ 
%\widetilde u(t)=\left\{\begin{array}{ll}
  %   \widetilde u_0, & t<t_0,\\
  %   \widetilde u_1, & t\in[t_0,t_1),\\
  %   \widetilde u_2, & t\ge t_1,
  %   \end{array}\right.\ \widetilde u_0=u_0,\ \widetilde u_2=u_1,
%\end{equation}
%which verify 
\begin{equation}\label{ver}
\eps\in(0,\eps_0),\ \|x_{u_0}\|\le d,\ \|x_{u_1}\|\le d, \  \|x_{u_0}-x_{v}\|\ge r, \ \|x_{v}-x_{u_1}\|\ge r.
\end{equation}
\end{theorem}

\vskip0.2cm

\noindent {\bf Proof.} By formula (\ref{Ti}) one has
\begin{eqnarray}
   T_{u_0,u_1}^\eps- T_{u_0,v}^\eps-T_{v,u_1}^\eps &=& -\dfrac{1}{k}\ln\dfrac{\eps}{\beta\left(\|x_{u_0}-x_{u_1}\|+\alpha^{-1}(\eps)\right)}+\nonumber\\
&&
+\dfrac{1}{k}\ln\dfrac{\eps}{\beta\left(\|x_{u_0}-x_{v}\|+\alpha^{-1}(\eps)\right)}+\nonumber\\
&&+\dfrac{1}{k}\ln\dfrac{\eps}{\beta\left(\|x_{v}-x_{u_1}\|+\alpha^{-1}(\eps)\right)}=\nonumber\\
%&=&  -\dfrac{1}{k_{u_1}}\ln\dfrac{\eps \cdot \beta_{\widetilde u_2}\left(\|x_{\widetilde u_1}-x_{\widetilde u_2}\|+\alpha_{\widetilde u_1}^{-1}(\eps)\right)}{\beta_{ u_1}\left(\|x_{u_0}-x_{u_1}\|+\alpha_{u_0}^{-1}(\eps)\right)\cdot\eps}+\\
%& & + \dfrac{1}{k_{\widetilde  u_1}}\ln\dfrac{\eps}{\beta_{\widetilde u_1}\left(\|x_{ \widetilde u_0}-x_{ \widetilde  u_1}\|+\alpha_{\widetilde u_0}^{-1}(\eps)\right)}=\\
&=& -\ln  \dfrac{K}{\eps^{1/k}},\label{f4}
\end{eqnarray} 
where
$$
  K={\left(\dfrac{\beta\left(\|x_{ v}-x_{  u_1}\|+\alpha^{-1}(\eps)\right)}{\beta\left(\|x_{  u_0}-x_{   u_1}\|+\alpha^{-1}(\eps)\right)}\right)^{{1}/{k}}} 
{\left({\beta\left(\|x_{ u_0}-x_{  v}\|+\alpha^{-1}(\eps)\right)}\right)^{{1}/{k}}}.
$$
Observe that there exists $K_0>0$ such that $K\ge K_0$ for any functions $u_0,u_1,v$ that verify  (\ref{ver}) as long as $d>0$ and $r>0$ stay fixed. Therefore, it is possible to choose $\eps_0>0$ (which depends on just $d>0$ and $r>0$) to satisfy  ${K}/{\eps^{1/k}}>1$ for all $\eps\in(0,\eps_0)$. The proof is complete.\qed

\begin{figure}[h]\center
\includegraphics[scale=0.44]{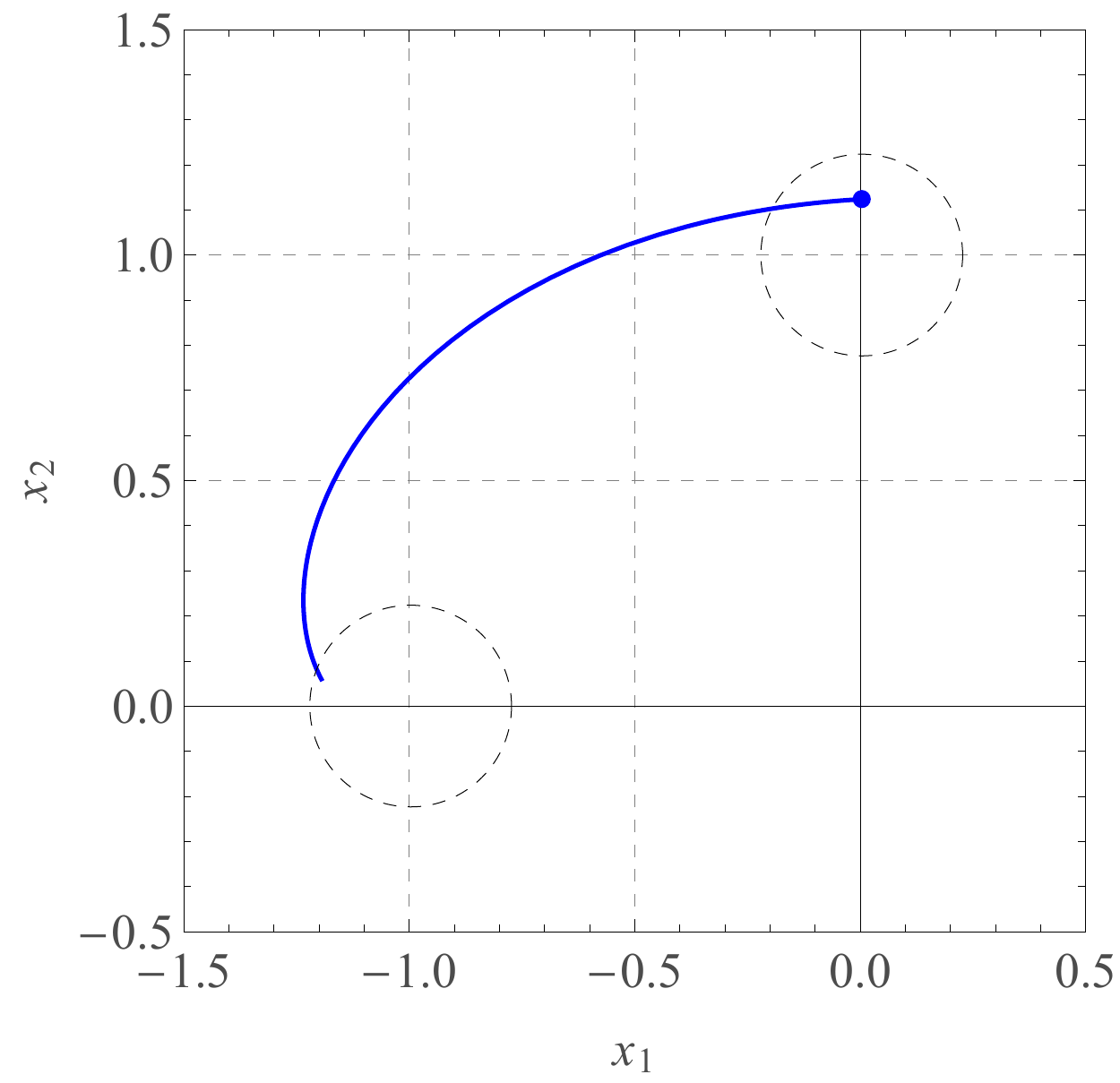}\ \ \ \
\includegraphics[scale=0.44]{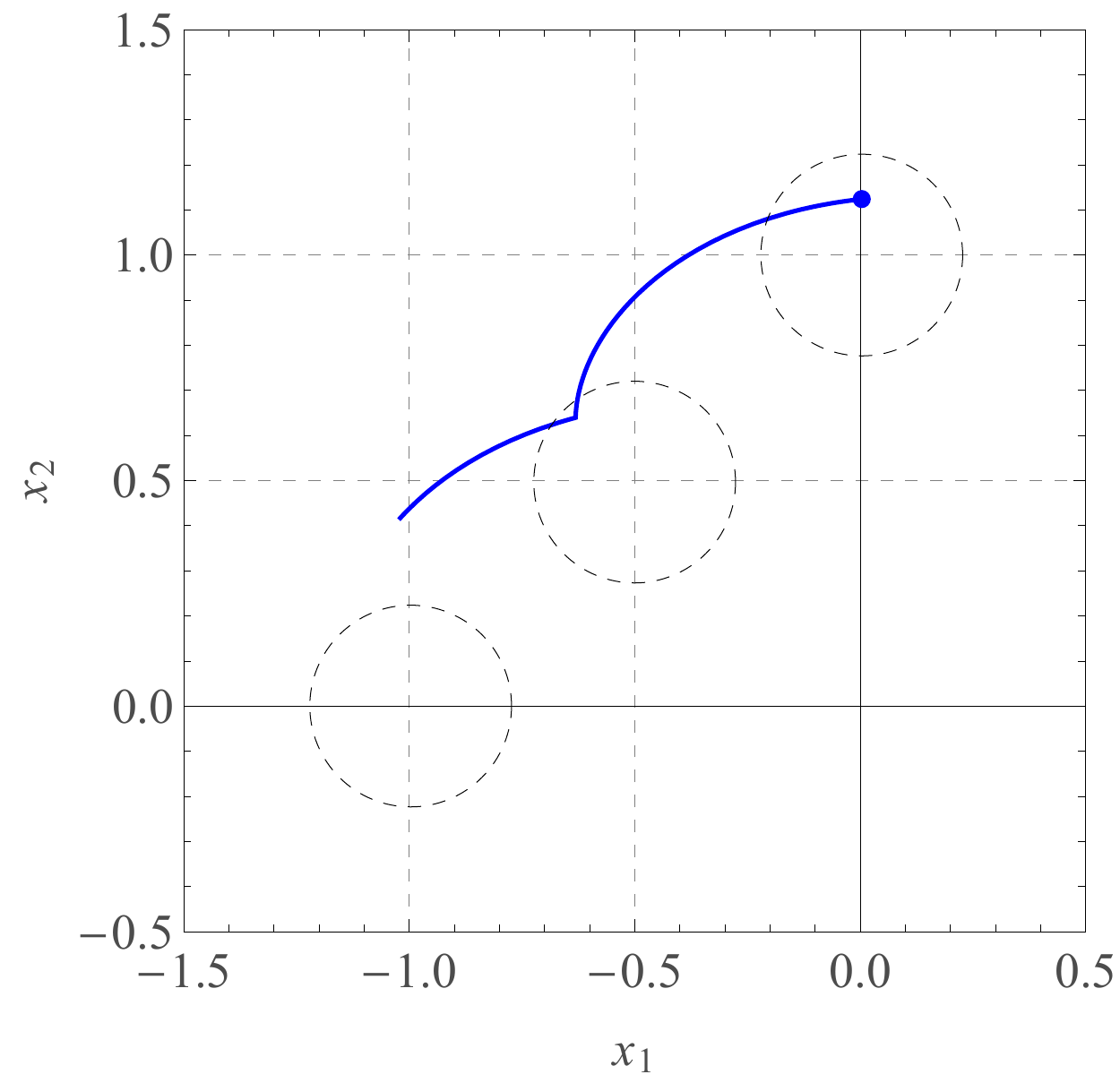}
\vskip-0.4cm
\caption{\footnotesize Solutions of switched system (\ref{ss}) with the initial condition in $N^{0.05}_{u_1}$ for the control inputs $\widetilde u(t)$ (Left) and $u(t)$ (Right) of Example~\ref{ex2}.} \label{fig2}
\end{figure}

\begin{ex}\label{ex2} In order to illustrate Theorem~\ref{thm41}, we refer to Example~\ref{ex1} again. Figure~\ref{fig2} shows the graphs of the solutions $x$ of (\ref{ss_ex1}) for two control inputs
$$
   \widetilde u(t)=\left\{\begin{array}{l}u_1,\quad t\in[0,T^\eps_{u_1,u_3}),\\
u_3,\quad t\ge T^\eps_{u_1,u_3},\end{array}\right.\ \ u(t)=\left\{\begin{array}{l}u_1,\quad t\in[0,T_{u_1,u_2}^\eps),\\
u_2,\quad t\in[T_{u_1,u_2}^\eps,T_{u_1,u_2}^\eps+T_{u_2,u_3}^\eps),\\
u_3,\quad t\ge T_{u_1,u_2}^\eps+T_{u_2,u_3}^\eps,\end{array}\right.
$$
over the time interval $[0,T^\eps_{u_1,u_3}]$. The plotting documents that $T_{u_1,u_2}^\eps+T_{u_2,u_3}^\eps$ turns out to be a longer time compared to $T^\eps_{u_1,u_3}.$
\end{ex}

\section{Conclusion} In this paper we considered a switched system of differential equations under the assumption that the time between two successive switchings is greater than a certain number $T$ called dwell time. We proved (Theorem~\ref{thm1}) that a suitable choice of the dwell time makes the solution stay within a required neighborhood $A_\eps$ of a so-called ideal attractor. We further proved that the solutions reach $A_\eps$  asymptotically if the initial conditions don't belong to $A_\eps.$ By doing that we obtained a new integral condition (Theorem~\ref{thm2}) for global stability which didn't seem to appear in the literature before. Finally, we addressed a case study where the Lyapunov functions of different subsystems are just shifts of one another. Here we used the dwell time formulas from Theorem~\ref{thm1} to estimate the time that the trajectories need to go from the neighborhood of an equilibrium of one subsystem to the neighborhood of an equilibrium of another subsystem (i.e. we considered a switched system with two discrete states). We proved (Theorem~\ref{thm41}) that adding more discrete states makes this travel time longer. Examples~\ref{ex1} and \ref{ex2} show that our theoretical conclusions agree with numeric simulations.

\section{Acknowledgements}

\noindent The first author is partially supported by NSF Grant CMMI-1436856.

\footnotesize
%\section{References}

% cut -----------
%\noindent References should be listed in alphabetical order.
%Author's initials should precede their names. References are
%indicated in the text by arabic numbers enclosed in square
%brackets. All the references should be cited in the paper.
%\medskip
% cut -----------

\noindent
email:\ journal@monotone.uwaterloo.ca\\
http://monotone.uwaterloo.ca/$\sim$journal/

\end{document}